\documentclass{article}
\usepackage{amsfonts,amsmath,amssymb}
\usepackage{graphics, epsfig}
\usepackage{color}
\allowdisplaybreaks
\newcommand{\dashiint}{{\dashint\!\!\dashint}}
\begin{document}


\newtheorem{theorem}{Theorem}[section]
\newtheorem{proposition}{Proposition}[section]
\newtheorem{lemma}{Lemma}[section]
\newtheorem{corollary}{Corollary}[section]
\newtheorem{remark}{Remark}[section]

\renewcommand{\thesection}{\arabic{section}}
\renewcommand{\theequation}{\thesection.\arabic{equation}}
\renewcommand{\thetheorem}{\thesection.\arabic{theorem}}
\numberwithin{equation}{section}
\numberwithin{theorem}{section}
\numberwithin{proposition}{section}
\numberwithin{lemma}{section}
\numberwithin{remark}{section}
\setcounter{secnumdepth}{3}



\newcommand{\cl}{\centerline}
\newcommand{\sms}{\smallskip}
\newcommand{\ms}{\medskip}
\newcommand{\bs}{\bigskip}
\newcommand{\noi}{\noindent}
\newcommand{\itl}[1]{\textit{#1}}
\newcommand{\blf}[1]{\textbf{#1}}
\newcommand{\dsty}{\displaystyle}
\newcommand{\txty}{\textstyle}
\newcommand{\ssty}{\scriptstyle}
\newcommand{\tty}{\texttt}


\newcommand\Par{\mathhexbox278\,}


\newcommand{\al}{\alpha}
\newcommand{\Al}{\Alpha}
\newcommand{\be}{\beta}
\newcommand{\Be}{\Beta}
\newcommand{\Gm}{\Gamma}
\newcommand{\gm}{\gamma}
\newcommand{\dl}{\delta}
\newcommand{\Dl}{\Delta}
\newcommand{\lm}{\lambda}
\newcommand{\Lm}{\Lambda}
\newcommand{\kp}{\kappa}
\newcommand{\varep}{\varepsilon}
\newcommand{\eps}{\epsilon}
\newcommand{\vp}{\varphi}
\newcommand{\sig}{\sigma}
\newcommand{\Sig}{\Sigma}
\newcommand{\om}{\omega}
\newcommand{\Om}{\Omega}
\newcommand{\uom}{\mbox{\boldmath$\omega$}}
\newcommand{\btau}{\mbox{\boldmath$\tau$}}
\newcommand{\bnu}{\mbox{\boldmath$\nu$}}
\newcommand{\up}{\upsilon}
\newcommand{\z}{\zeta}


\newcommand{\df}[1]{\buildrel\mbox{\small def}\over{#1}}
\newcommand{\op}[1]{\buildrel\mbox{\tiny o}\over{#1}}
\newcommand{\db}{\prime\prime}
\newcommand{\bsl}{\backslash}
\newcommand{\lb}{\lbrack\!\lbrack}
\newcommand{\rb}{\rbrack\!\rbrack}
\newcommand\la{\langle}
\newcommand\ra{\rangle}
\newcommand{\ev}{\equiv}
\newcommand{\nev}{\not\equiv}
\newcommand{\nn}{\mathbb{N}}
\newcommand{\qq}{\mathbb{Q}}
\newcommand{\zz}{\mathbb{Z}}
\newcommand{\rr}{\mathbb{R}}
\newcommand{\rn}{\rr^N}
\newcommand{\cc}{\mathbb{C}}
\newcommand{\id}{\mathbb{I}}
\newcommand{\bo}{\mathbb{O}}

\newcommand{\amsb}[1]{\mathbb{#1}}
\newcommand{\mcl}[1]{\mathcal{#1}}
\newcommand{\bl}[1]{\mathbf{#1}}
\newcommand{\ov}[1]{\overline{#1}}
\newcommand{\wt}[1]{\widetilde{#1}}
\newcommand{\wh}[1]{\widehat{#1}}

\newcommand{\llra}{\leftrightarrow}
\newcommand{\lra}{\longrightarrow}
\newcommand{\LLR}{\Longleftrightarrow}
\newcommand{\LRA}{\Longrightarrow}
\newcommand{\LLA}{\Longleftarrow}


\newcommand{\bbox}{\vrule height.6em width.6em 
depth0em} 
\newcommand{\os}{\vbox{\hrule \hbox{\vrule 
height.6em depth0pt 
\hskip.6em \vrule height.6em depth0em}
\hrule}} 


\newcommand{\dvg}{\operatorname{div}}
\newcommand{\curl}{\operatorname{curl}}
\newcommand{\supp}{\operatorname{supp}}
\newcommand{\essup}{\operatornamewithlimits{ess\,sup}}
\newcommand{\essinf}{\operatornamewithlimits{ess\,inf}}
\newcommand{\essosc}{\operatornamewithlimits{ess\,osc}}
\newcommand{\osc}{\operatornamewithlimits{osc}}
\newcommand{\sign}{\operatorname{sign}}
\newcommand{\loc}{\operatorname{loc}}
\newcommand{\diam}{\operatorname{diam}}
\newcommand{\dist}{\operatorname{dist}}
\newcommand{\card}{\operatorname{card}}
\newcommand{\meas}{\operatorname{meas}}
\newcommand{\spn}{\operatorname{span}}
\newcommand{\dtm}{\operatorname{det}}
%


\newcommand{\overlim}{\mathop{\overline{\lim}}\limits}
\newcommand{\underlim}{\mathop{\underline{\lim}}\limits}
\newcommand{\ttop}[2]{\genfrac{}{}{0pt}{}{#1}{#2}}
\newcommand{\bcu}{\mathop{\txty{\bigcup}}\limits}
\newcommand{\bca}{\mathop{\txty{\bigcap}}\limits}
\newcommand{\bsu}{\mathop{\txty{\sum}}\limits}
\newcommand{\pro}{\mathop{\txty{\prod}}\limits}


\newcommand{\pl}{\partial}
\newcommand{\ptt}{\frac{\pl}{\pl t}}
\newcommand{\ppx}{\frac\pl{\pl x}}
\newcommand{\dds}{\frac d{ds}}
\newcommand{\ddt}{\frac d{dt}}

\newcommand{\intl}{\int\limits}
\newcommand{\iintl}{\iint\limits}
\def\Xint#1{\mathchoice
    {\XXint\displaystyle\textstyle{#1}}%
    {\XXint\textstyle\scriptstyle{#1}}%
    {\XXint\scriptstyle\scriptscriptstyle{#1}}%
    {\XXint\scriptscriptstyle\scriptscriptstyle{#1}}%
    \!\int}
\def\XXint#1#2#3{\setbox0=\hbox{$#1{#2#3}{\int}$}
    \vcenter{\hbox{$#2#3$}}\kern-0.5\wd0}
\def\bint{\Xint-}
\def\dashint{\Xint{\raise4pt\hbox to7pt{\hrulefill}}}


\newcommand{\ovl}[3]{\int_{#1}^{#2}\kern-#3pt\raise4pt\hbox to7pt{\hrulefill}\ }

\newcommand{\ovll}[3]{\intl_{#1}^{#2}\kern-#3pt\raise4pt\hbox to7pt{\hrulefill}\ }

\newcommand{\tvl}[2]{\iint_{#1}\kern-#2pt\raise4pt\hbox to7pt{\hrulefill}\ }



\newcommand{\omt}{\Om_T}
\newcommand{\plo}{\partial\Omega}
\newcommand{\ovo}{\bar{\Om} }

%
\newcommand{\ci}[1]{C^\infty\!\left({#1}\right)}
\newcommand{\cio}[1]{C_o^\infty\!\left({#1}\right)}
\newcommand{\lloc}[1]{L_{\loc}\!\left({#1}\right)}
\newcommand{\xy}{|x-y|}


\newcommand{\intom}{\intl_{\Om}}
\newcommand{\intbo}{\intl_{\plo}}
\newcommand{\inom}{\int_{\Om}}
\newcommand{\inbo}{\int_{\plo}}
\newcommand{\intrn}{\intl_{\rn}}


\newcommand{\bye}{\end{document}}



%
%
%

%
%
%
%
%

%
%

\input harnack_mono.mac
\title{A sufficient condition for the continuity of
solutions to a logarithmic diffusion equation}
\author{Naian Liao\footnote{Supported by Chongqing University Grant No. 106112015CDJXY100006}}
\date{}
\maketitle
\begin{abstract}
This note gives a first sufficient condition that insures
a non-negative, locally bounded, local solution to a logarithmically singular parabolic equation
is continuous at a vanishing point and an estimate of the modulus of 
continuity is given. Moreover, an estimate of
the Hausdorff measure of the set of discontinuity is established.
\end{abstract}

\section{Introduction and Main Results}
Let $E$ be an open set in $\rn$. For $T>0$, let 
$E_T$ denote the cylindrical domain $E\times(0,T]$. 
Consider the quasi-linear, parabolic 
differential equation
\begin{equation}\label{Eq:1:1}
u_t-\Dl\ln u = 0\quad
\text{ weakly in }\> E_T.
\end{equation}
This equation is singular since its modulus of 
ellipticity $u^{-1}\to\infty$ as $u\to 0$. 
A non-negative function $u$ satisfying
\begin{equation*}
u\in C_{loc}(0,T;L^2_{loc}(E)),\quad\ln u\in 
L^2_{loc}(0,T;W^{1,2}_{loc}(E))
\end{equation*}
is called a local, weak sub(super)-solution to (\ref{Eq:1:1}) 
if for every compact set $K\subset E$ and every 
sub-interval $[t_1,t_2]\subset(0,T]$
\begin{equation*}
\int_K u\vp dx\Big|_{t_1}^{t_2}+\int_{t_1}^{t_2}
\int_K\Big(-u\vp_t+\frac{Du}uD\vp\Big)dxdt\le\,(\ge)\,0
\end{equation*}
for all non-negative testing functions
\begin{equation*}
\vp\in W^{1,2}_{loc}(0,T;L^2(K))\cap 
L^2_{loc}(0,T;W^{1,2}_o(K)).
\end{equation*}
A function $u$ that is both a local, weak sub-solution and 
a local, weak super-solution is a local, weak solution.

For $\rho>0$ we denote by $K_{\rho}(y)$ the cube centered at $y$
with side length $\rho$. If $y=0$ we use $K_\rho$. 
For $\theta>0$ introduce the cylinder with ``vertex" at $(0,0)$
\[Q_{\rho}(\theta)=K_{\rho}\times(-\theta\rho^2,0].\]
If $\theta=1$ we use $Q_\rho$.
Also a cylinder with ``vertex" at $(y,s)$ is
\[(y,s)+Q_{\rho}(\theta)=K_{\rho}(y)\times(s-\theta\rho^2,s].\]
Assume $u$ is a locally bounded, local solution.
 Let us suppose $\rho>0$ is so small that the cylinder $(y,s)+Q_{\rho}\subset E_T$.
Up to a translation we may assume $(y,s)=(0,0)$ and let $$\om=\essosc_{Q_\rho}u.$$
Without loss of generality we assume $\om\le1$ such that
\[Q_{\rho}(\om)\subset Q_{\rho}\quad\text{and}\quad\essosc_{Q_{\rho}(\om)}u\le \om.\]
Suppose in addition to the notion of solution that
\begin{equation}\label{Eq:1:2}
D\ln u\in L^p_{loc}(E_T)\quad\text{for some }p>\frac{N+2}2.
\end{equation}
Note that when $N=1$ the integrability condition \eqref{Eq:1:2} is inherent in the
notion of solution while in other cases it has to be imposed. Accordingly we define the quantity 
\[
I_{p,\rho}(y,s)=\rho\bigg(\dashiint_{(y,s)+Q_{\rho}}|D\ln u|^{p}\,dxdt\bigg)^{\frac1p}
\]
and $I_{p,\rho}=I_{p,\rho}(0,0)$.
Then we have the following main theorem.
\begin{theorem}\label{Thm:1}
Let $u$ be a non-negative, locally bounded, local solution to \eqref{Eq:1:1} 
and assume \eqref{Eq:1:2} is satisfied. 
Then there exist constants $\bar{C}>1$ and $\al\in(0,1)$ depending only
on $N$, such that for any $\mu\in(0,1)$ and $0<r<\rho\le R_o$ we have 
\[
\essosc_{Q_{r}(\om)}u\le\bar{C}\bigg[\om\bigg(\frac{r}{R_o}\bigg)^{(1-\mu)\al}+I_{p,R^{1-\mu}_o r^\mu}\bigg]
\]
In particular, the solution $u$ is continuous at the origin provided
\begin{equation}\label{Eq:1:3}
\limsup_{r\to 0}I_{p,r}=0.
\end{equation}
\end{theorem}
\begin{remark}
Strictly speaking, we need the convention that the function $\rho\to I_{p,\rho}$
is non-decreasing. In order to validate that, we need only to take
\[
\tilde{I}_{p,\rho}=\sup_{0<\tau<\rho}I_{p,\tau}
\]
in Theorem \ref{Thm:1}
\end{remark}
\begin{remark}
For $\lm\ge0$, $T>0$ and $N\ge3$ the explicit solution
\begin{equation}\label{Eq:ex_sol}
u(x,t)=\frac{2(N-2)(T-t)^{\frac{N}{N-2}}}{\lm+(T-t)^{\frac{2}{N-2}}|x|^2}
\end{equation}
is continuous up to its extinction time $T$.
One verifies that when $\lm>0$ and for any fixed $x_o$,
 there is a positive constant $C(x_o, \lm, N, p)$ such that
\[
I_{p,r}(x_o,T)\le C r^{\frac{4}{(N-2)p}+2}\to0\quad\text{as }r\to0.
\]
When $\lm=0$, it gives an unbounded solution
which, in particular, is discontinuous at $x=0$.
Condition \eqref{Eq:1:3} is verified everywhere except for $x=0$.
A direct calculation shows that
\[D\ln u\in L^p_{loc}(\rr^N\times\rr_+)\quad\text{for any}\quad\frac{N+2}2<p<N.\]
 Furthermore, there exists some positive constant $C(N, p)$ such that for every $t<T$
\begin{equation*}
I_{p,r}(t,0)=\left\{
\begin{array}{ll}
\dsty C(N,p),\quad \frac{N+2}2<p<N;\\
\dsty\infty, \quad p\ge N.
\end{array}\right.
\end{equation*}
Hence the condition \eqref{Eq:1:2} alone is not sufficient to ensure continuity.
\end{remark}
Now define the set $S\subset E_T$ to consist of all discontinuous points of a local solution $u$ and 
\[
S_o=\left\{(y,s)\in E_T: \limsup_{\rho\to0}\frac{1}{\rho^{N+2-p}}\iint_{(y,s)+Q_\rho}|D\ln u|^p\,dxdt>0\right\}.
\]
As a direct consequence of Theorem \ref{Thm:1} it is straightforward to see that $S\subset S_o$.
Moreover we are going to obtain an estimate of the Hausdorff measure 
of the set $S_o$.

The {\it parabolic Hausdorff measure} $P_k$ is defined in a way similar to the usual Hausdorff measure $H_k$
but using the parabolic metric on $\rr^N\times \rr$. For any set $U\subset \rr^N\times\rr$
and $k\ge0$ we define 
\[
P_k(U)=\lim_{\dl\to0}P_k^{\dl}(U),
\]
where
\[
P_k^{\dl}(U)=\inf\left\{\sum_{i=1}^{\infty} r_i^k: \, U\subset\bigcup_{i}[(y_i,t_i)+Q_{r_i}],\, r_i<\dl\right\}.
\]
$P_k$ so defined is an outer measure whose $\sig$-algebra 
contains all Borel sets of $\rr^N\times\rr$ (Chapter 2, \cite{EG}).
It should be pointed out that the parabolic Hausdorff measure dominates 
the usual Hausdorff measure in the sense that there is some constant $C(N,\, k)$ 
such that for any subset $U$ of $\rr^N\times\rr$ one has
\[
H_k(U)\le CP_k(U).
\]

Regarding the Hausdorff measure of the discontinuity set we have the following 
consequence of Theorem \ref{Thm:1}.
\begin{theorem}\label{Thm:2}
Let $u$ be a non-negative, locally bounded, local solution to \eqref{Eq:1:1} 
and assume \eqref{Eq:1:2} is satisfied. Then we have
\begin{equation*}
\begin{array}{ll}
P_{N+2-p}(S_o)=0, \quad N>1\quad\text{and}\\
P_1(S_o)=0,\quad N=1.
\end{array}
\end{equation*}
\end{theorem}
\begin{remark}
When $N=1$ the possible discontinuous points
of a non-negative, locally bounded, local solution to \eqref{Eq:1:1} cannot occupy a line in $\rr^2$. Generally
one gets less discontinuity as the $L^p$ integrability of $D\ln u$ increases
and eventually, the solution is continuous at every point if one has
$p\ge N+2$.
\end{remark}
\subsection{Novelty and Significance}
Equation \eqref{Eq:1:1} describes the evolution of the Ricci flow for complete $\rr^2$ (\cite{W}).
It also arises from modeling the thickness of a viscous liquid thin film that lies on 
a rigid plate under the influence of the van der Waals force (\cite{WD}).

Physical and geometric motivations of \eqref{Eq:1:1} make sense mainly for $N=2$,
but the problem is intriguing in the effort to shed light on the structural properties of 
singular diffusion equations.

Questions concerning both existence and non-existence of solutions to the Cauchy
problem of \eqref{Eq:1:1} and its related elliptic equation are investigated
in \cite{DDP, DDP1, DD, DD1, Hsu1, Hui1,V2} (just mention few).

The study of local behavior of local solutions to \eqref{Eq:1:1} has been 
initiated in \cite{DBGL1,DBGL2}. Equation \eqref{Eq:1:1} can be viewed
as a formal limit of the porous medium equation
\[
u_t-\dvg (u^{m-1}Du)=0\quad\text{ as }m\to0.
\]
A proof of H\"older continuity for non-negative, locally bounded, local solutions
to the porous medium equation can be found in Appendix~B of \cite{DGV}.
However, the local behavior of local solutions to \eqref{Eq:1:1} presents many striking differences
from that of local solutions to the porous medium equation.
See \cite{L} for more detailed discussion.

It was shown in \cite{DBGL1} that if one assumes that
\[
u\in L^r_{loc}(E_T)\quad\text{for some}\quad r>\max\{1,\frac{N}2\},
\]
then $u$ is locally bounded.
If in addition one assumes that
\[
\ln u\in L_{loc}^{\infty}(0,T;L^p_{loc}(E))\quad\text{for some}\quad p>N+2,
\]
then a Harnack-type inequality is established and thus, if the solution
does not vanish identically on a hyperplane normal to the time axis,
then the equation \eqref{Eq:1:1} is neither degenerate nor singular in 
a backward cylinder with its vertex on the hyperplane. 
As a result $u$ is a classical solution in such a cylinder. In fact, it is shown in \cite{DBGL3}
that under such circumstances the solution is analytic in space variables while infinitely differentiable in time.

Nevertheless, these results do not explain why some explicit solutions, \eqref{Eq:ex_sol} for example, could be 
continuous up to their extinction time. Theorem \ref{Thm:1} gives a first sufficient
condition that insures continuity at a vanishing point of $u$, and 
an explicit estimate of the modulus of continuity is given.
Moreover,
we establish in Theorem \ref{Thm:2} an estimate on the Hausdorff measure of the set of 
discontinuity of $u$. 

Those effort being made, it is interesting to ask whether the higher integrability conditon \eqref{Eq:1:2} of $D\ln u$
for $N>1$ can be obtained from the notion of solution and whether
the condition \eqref{Eq:1:3} is necessary for a point to be a continuity point of $u$.
Last but not least, can we construct an explicit bounded solution with discontinuity? 
In \cite{V1} when $N=1$, a solution discontinuous on a line segment
was constructed.  However, the notion of solution used seems different
from this note, since our results indicate such a phenomenon  is not allowed
for our solutions. 
\\
\\
{\it Acknowledgement.} This paper was finalized during my visit to Vanderbilt University
in November 2016. I am grateful for many helpful discussions with Professor Emmanuele DiBenedetto
and Professor Ugo Gianazza. Professor Gianazza also read carefully the early version of this paper
and came up with a lot of valuable comments. I am really indebted to both of them.
\section{Proof of Theorem \ref{Thm:2} Assuming Theorem \ref{Thm:1}}
The proof of Theorem \ref{Thm:2} is based on the following
\begin{proposition}\label{Prop:1}
Let $f\in L^1_{loc}(\rr^{N+1})$, suppose $0\le s<N+2$ and define
\[
\Lm_s=\left\{(y,s)\in \rr^{N+1}: \limsup_{\rho\to0}\frac{1}{\rho^s}\iint_{(y,s)+Q_\rho}|f|\,dxdt>0\right\}.
\]
Then
\[
P_s(\Lm_s)=0.
\]
\end{proposition}

This is a parabolic counterpart of a similar result shown in \cite{EG} (Theorem 3, p.77). 
Now we are ready to present\\
{\bf Proof of Theorem \ref{Thm:2}.} When $N>1$, since we assume
$$D\ln u\in L^p_{loc}(E_T)\quad\text{for some }p>\frac{N+2}{2},$$
a straightforward application of Proposition \ref{Prop:1}
yields the desired conclusion.

When $N=1$, the notion of solution gives $$D\ln u\in L^2_{loc}(E_T)$$ and 
by the H\"older inequality with $p<2$
\[
\bigg[\frac{1}{\rho^{3-p}}\iint_{Q_\rho}|D\ln u|^{p}\,dxdt\bigg]^{\frac1p}
\le\bigg[\frac{1}{\rho}\iint_{Q_\rho}|D\ln u|^2\,dxdt\bigg]^{\frac12}.
\]
Thus
\[
S_o\subset \left\{(y,s)\in E_T: \limsup_{\rho\to0}\frac{1}{\rho}\iint_{(y,s)+Q_\rho}|D\ln u|^2\,dxdt>0\right\}
\]
and again by Proposition \ref{Prop:1} we obtain
\[P_1(S_o)=0.\]
This finishes the proof.\hfill\bbox

\noi The rest of the note is devoted to proving Theorem \ref{Thm:1}.
\section{Some Preliminary Estimates}
\subsection{Energy Estimates}
\begin{proposition}
Let $u$ be a local, weak super-solution to \eqref{Eq:1:1}. Then there is a
positive constant $\gm$ depending only on $N$ such that for every
cylinder $(y,s)+Q_{\rho}(\theta)\subset E_T$, every $k\in \rr_+$, and every non-negative, piecewise smooth
cutoff function $\z$ vanishing on $\pl K_{\rho}(y)$,
\begin{equation*}
\begin{aligned}
\essup_{s-\theta\rho^2<t<s}&\frac12\int_{K_\rho(y)}(u-k)_-^2\z^2\,dx
+k^{-1}\iint_{(y,s)+Q_{\rho}(\theta)}|D[(u-k)_-\z]|^2\,dxdt\\
&\le \frac12\int_{K_{\rho}(y)}(u-k)_-^2\z^2(x,s-\theta\rho^2)\,dx\\
&+\iint_{(y,s)+Q_{\rho}(\theta)}(u-k)_-^2 \z|\z_t|\,dxdt\\
&+k^{-1}\iint_{(y,s)+Q_{\rho}(\theta)}(u-k)_-^2 |D\z|^2\,dxdt\\
&+2\iint_{(y,s)+Q_{\rho}(\theta)}|D\ln u|(u-k)_- |D\z|\z\,dxdt.
\end{aligned}
\end{equation*}
\end{proposition}
{\bf Proof.} We may assume $(y,s)=(0,0)$. In the weak formulation for super-solutions
to $\eqref{Eq:1:1}$, we take the test function
\[
\vp=-(u-k)_-\z^2
\]
over the cylinder
\begin{equation*}
Q_t=K_\rho\times(-\theta\rho^2,t] \quad\text{ for }\quad 
t\in(-\theta\rho^2,0],
\end{equation*} 
modulo a standard Steklov averaging process.
This gives
\begin{align*}
&-\iint_{Q_t}u_\tau (u-k)_-\z^2\,dxd\tau-\iint_{Q_t}\z^2 D\ln u D(u-k)_-\,dxd\tau\\
&=\iint_{Q_t}2\z (u-k)_- D\ln u D\z\,dxd\tau.
\end{align*}
The first term on the left-hand side is estimated by
\begin{align*}
&-\iint_{Q_t}u_\tau (u-k)_-\z^2\,dxd\tau\\
&\ge\frac12\int_{K_\rho}(u-k)_-^2\z^2(x,t)\,dx
-\frac12\int_{K_\rho}(u-k)_-^2\z^2(x,-\theta\rho^2)\,dx\\
&\quad-\iint_{Q_{\rho}(\theta)}(u-k)_-^2\z|\z_\tau|\,dxd\tau,
\end{align*}
while the second term is estimated by
\begin{align*}
-\iint_{Q_t}\z^2 D\ln u D(u-k)_-\,dxd\tau
\ge k^{-1}\iint_{Q_{\rho}(\theta)}|D(u-k)_-|^2\z^2\,dxd\tau.
\end{align*}
Next the term on the right side is 
\[
\iint_{Q_t}2\z (u-k)_- D\ln u  D\z\,dxd\tau\le 2\iint_{Q_{\rho}(\theta)}|D\ln u|(u-k)_- |D\z|\z\,dxd\tau
\]
Combining all these estimates yields the conclusion.\hfill\bbox
\begin{proposition}\label{Prop:2:3}
Let $u$ be a local, weak sub--solution to 
(\ref{Eq:1:1}) in $E_T$.  
There exists a positive constant $\gm=\gm(N)$, 
such that for every cylinder
$(y,s)+Q_{\rho}(\theta)\subset E_T$, every 
$k\in\rr_+$, and every non-negative, piecewise smooth 
cutoff function $\z$ vanishing on $\pl K_\rho (y)$, 
\begin{equation}\label{Eq:2:13}
\begin{aligned}
\essup_{s-\theta\rho^2<t\le s}&\int_{K_\rho(y)}\ukp^2\z^2(x,t)dx\\
&\quad -\int_{K_\rho(y)}\ukp^2
\z^2(x,s-\theta\rho^2)dx\\
&\quad + \iint_{(y,s)+Q_\rho(\theta)}\frac{|D[(u-k)_+\z]|^2}{u}dxdt\\
&\le \gm\iint_{(y,s)+Q_\rho(\theta)}\ukp^2
\z|\z_t|dxdt\\
&\quad+\gm\iint_{(y,s)+Q_\rho(\theta)}\frac{\ukp^2}{u}|D\z|^2dxdt.
\end{aligned}
\end{equation}
\end{proposition}
\noi{\bf Proof.}
After a translation may assume $(y,s)=(0,0)$. 
Take the test function $\vp=\ukp\z^2$ over $Q_t$ 
modulo a standard Steklov averaging process, and perform 
standard calculations. The various integrals are extended over the 
set $[u>k]$ and since $k>0$, they are all well defined. \hfill\bbox
\subsection{A Logarithmic Estimate for Sub-Solutions}
Introduce the logarithmic function
\begin{equation}\label{Eq:B:7:1}
\psi(u)=\ln^+\Big[\frac{H}{H-\ukp +c}\Big]
\end{equation}
where 
\begin{equation*}
H=\essup_{(y,s)+Q_{\rho}(\theta)}\ukp,\quad  0<c< \min\{1;H\},
\end{equation*} 
and for $s>0$
\begin{equation*}
\ln^+ s=\max\{\ln s;0\}.
\end{equation*}
In the cylinder $(y,s)+Q_{\rho}(\theta)$ take a non-negative, 
piecewise smooth cutoff function $\z$  independent of $t$.
\begin{proposition}\label{Prop:B:7:1} 
Let $u$ be a non-negative, locally bounded, local, weak 
sub-solution to equation 
(\ref{Eq:1:1}) in $E_T$. There exists a  constant $\gm$, 
depending only on the $N$, such that for every
cylinder 
\begin{equation*}
(y,s)+Q_{\rho}(\theta)\subset E_T 
\end{equation*}
and for every level $k\ge0$ we have
\begin{align}
&\sup_{s-\theta\rho^2<t<s} \int_{K_\rho(y)}
\psi^2(u)(x,t)\z^2(x)dx\nonumber\\
&\le\int_{K_\rho(y)}\psi^2(u)(x,s-\theta\rho^2)\z^2(x)dx
+\gm \iint_{(y,s)+Q_\rho(\theta)}\frac{\psi(u)}u |D\z\big|^2dxdt.
\label{Eq:B:7:2}
\end{align}
\end{proposition}
\noi{\bf Proof.}  
Take $(y,s)=(0,0)$ and work within the cylinder $Q_t$ 
introduced before in the energy estimates. In the weak formulation 
of (\ref{Eq:1:1}) take the testing function
\begin{equation*}
\vp=\frac{\pl}{\pl u}\big[\psi^2(u)\big]\z^2
= 2\psi\psi^\prime\z^2.
\end{equation*}
By direct calculation
\begin{equation*}
\big[\psi^2(u)\big]^{\db}=2(1+\psi)\psi^{\prime2}
\in L^\infty_{loc}(E_T)
\end{equation*}
which implies that such a $\vp$ is an admissible testing 
function, modulo a Steklov averaging process.
Since $\psi(u)$ vanishes on the set where $\ukp=0$
\begin{equation*}
\iint_{Q_t} u_\tau[\psi^2]^\prime\z^2dxdt
=\int_{K_\rho}\psi^2(x,t)\z^2dx -
\int_{K_\rho}\psi^2(x,-\theta\rho^2)\z^2dx.
\end{equation*}
As for the remaining term
\begin{align*}
&\iint_{Q_t}\frac{Du}u\cdot D\vp dxd\tau\\
&\ge 2\iint_{Q_t}(1+\psi)\psi^{\prime2} \frac{|Du|^2}u\z^2dxd\tau
-4\iint_{Q_t}\frac{|Du|}u\psi\psi^\prime\z|D\z|dxd\tau\\
&\ge \iint_{Q_t}(1+\psi)\psi^{\prime2} \frac{|Du|^2}u\z^2dxd\tau-\gm\iint_{Q_t}\frac{\psi}u|D\z|^2dxd\tau.
\end{align*}
Collecting these estimates establishes the proposition.\hfill\bbox
\section{DeGiorgi-type Lemmas}
For a cylinder $(y,s)+Q_{2\rho}(\theta)\subset E_T$ 
denote by $\mu_\pm$ and $\om$, numbers satisfying 
\begin{equation*}
\mu_+\ge\essup_{[(y,s)+Q_{2\rho}(\theta)]}u,\qquad
\mu_-\le\essinf_{[(y,s)+Q_{2\rho}(\theta)]}u,\qquad
\om=\mu_+-\mu_-.
\end{equation*}
Denote by $\xi$ and $a$ fixed numbers in $(0,1)$.
\begin{lemma}\label{Lm:4:1}
Let $u$ be a non-negative, locally bounded, local, 
weak super-solution to \eqref{Eq:1:1}. Then there is a constant
$\nu_-$ depending on the data and $\theta,\,\xi,\,\om,\,a,$ such that if
\[
|[u\le\mu^-+\xi\om]\cap[(y,s)+Q_{2\rho}(\theta)]|\le\nu_-|Q_{2\rho}(\theta)|,
\]
then either
\[\xi\om\le I_{p,\rho}(y,s),\]
or
\[
u\ge \mu^-+a\xi\om\quad\text{ a.e. in }(y,s)+Q_{\frac{\rho}2}(\theta).
\]
\end{lemma}
{\bf Proof.} We may take $(y,s)=(0,0)$. Set 
\[
 \rho_n=\rho+\frac{\rho}{2^n},\quad K_n=K_{\rho_n},\quad Q_n=K_n\times(-\theta\rho_n^2,0]
\]
Consider a non-negative, piecewise smooth cutoff function on 
$Q_n$ of the form $\z(x,t)=\z_1(x)\z_2(t)$, where
\begin{equation*}
\begin{array}{lc}
{\dsty 
\z_1=\left\{
\begin{array}{ll}
1\>&\text{ in }\> K_{n+1}\\
{}\\
0\>&\text{ in }\>\rn-K_n
\end{array}\right .}\quad
&{\dsty |D\z_1|\le \frac1{\rho_n-\rho_{n+1}}
=\frac{2^{n+1}}{\rho} }\\
{}\\
{\dsty 
\z_2=\left\{
\begin{array}{ll}
0\>&\text{ for }\> t<-\theta \rho_n^2\\
{}\\
1\>&\text{ for }\> t\ge-\theta\rho_{n+1}^2
\end{array}\right .}\quad
&{\dsty 0\le \z_{2,t}\le\frac1{\theta(\rho_n^2-\rho_{n+1}^2)}
\le\frac{2^{2(n+1)}}{\theta \rho^2}}
\end{array}
\end{equation*} 
Now apply the energy estimate to $(u-k_n)_-$ in the cylinder $Q_n$
with $$k_n=\mu^-+a\xi\om+\frac{1-a}{2^n}\xi\om$$
to obtain
\begin{align*}
\essup_{-\theta\rho^2<t<0}&\int_{K_n}(u-k_n)_-^2\z^2(x,t)\,dx+\frac1{\mu^-+\xi\om}\iint_{Q_n}|D[(u-k_n)_-\z]|^2\,dxdt\\
&\le\frac{4^n}{\theta\rho^2}\iint_{Q_n}(u-k_n)_-^2\,dxdt+\frac{4^n}{(\mu^-+a\xi\om)\rho^2}\iint_{Q_n}(u-k_n)_-^2\,dxdt\\
&+4^n\frac{M\xi\om}{\rho}|[u<k_n]\cap Q_n|^{1-\frac1p}
\end{align*}
where
\[
M=\bigg(\iint_{Q_\rho}|D\ln u|^{p}\,dxdt\bigg)^{\frac1p}.
\]
Let $A_n=[u<k_n]\cap Q_n$. By the standard parabolic embedding theorem (Proposition 3.1, Chapter 1 of \cite{DB}), 
we obtain
\begin{align*}
&\bigg(\frac{1-a}{2^{n+1}}\xi\om\bigg)^2|A_{n+1}|\\
&\le\iint_{Q_{n}}(u-k_n)_-^2\z^2\,dxdt\\
&\le \bigg(\iint_{Q_{n}}[(u-k_n)_-\z]^{2\frac{N+2}N}\,dxdt\bigg)^{\frac{N}{N+2}}|A_n|^{\frac2{N+2}}\\
&\le\bigg(\iint_{Q_n}|D[(u-k_n)_-\z]|^2\,dxdt\bigg)^{\frac{N}{N+2}}\\
&\times\bigg(\essup_{-\theta\rho^2<t<0}\int_{K_n}(u-k_n)_-^2\z^2(x,t)\,dx\bigg)^{\frac2{N+2}}|A_n|^{\frac2{N+2}}\\
&\le4^n|A_n|^{\frac2{N+2}}(\mu^-+\xi\om)^{\frac{N}{N+2}}\\
&\times\bigg[\bigg(\frac1{\theta}
+\frac1{\mu^-+a\xi\om}\bigg)\frac{(\xi\om)^2}{\rho^2}|A_n|
+\frac{M\xi\om}{\rho}|A_n|^{1-\frac1p}\bigg].
\end{align*}
Setting
\[
Y_n=\frac{|A_n|}{|Q_n|},
\]
 we have
\begin{align*}
Y_{n+1}&\le\frac{\gm 4^{2n}}{(1-a)^2}
\bigg[\bigg(\bigg(\frac{\mu^-+\xi\om}{\theta}\bigg)^{\frac{N}{N+2}}
+\gm_o\bigg(\frac{\theta}{\mu^-+a\xi\om}\bigg)^{\frac2{N+2}}\bigg)Y_{n}^{1+\frac2{N+2}}\\
&+\frac{I_{p,\rho}}{\xi\om}\bigg(\frac{\mu^-+\xi\om}{\theta}\bigg)^{\frac{N}{N+2}}Y_n^{1-\frac1p+\frac2{N+2}}\bigg]
\end{align*}
where
\[
\gm_o=\bigg(\frac{\mu^-+\xi\om}{\mu^-+a\xi\om}\bigg)^{\frac{N}{N+2}}.
\]
Suppose $I_{p,\rho}\le\xi\om$; we have
\begin{align*}
Y_{n+1}&\le\frac{\gm 4^{2n}}{(1-a)^2}
\bigg[\bigg(\frac{\mu^-+\xi\om}{\theta}\bigg)^{\frac{N}{N+2}}
+\gm_o\bigg(\frac{\theta}{\mu^-+a\xi\om}\bigg)^{\frac2{N+2}}\bigg]Y_{n}^{1+\be}
\end{align*}
where
\[
\be=\frac{2}{N+2}-\frac1p>0.
\]
It follows from Lemma 4.1 of Chapter 1 of \cite{DB} that $Y_n$ tend to $0$ provided
\begin{equation}\label{mu}
Y_o\le\nu_-\df{=}A^{-\frac{1}{\be}}16^{-\frac{1}{\be^2}},
\end{equation}
where 
\[
A=\frac{\gm}{(1-a)^2}\bigg[\bigg(\frac{\mu^-+\xi\om}{\theta}\bigg)^{\frac{N}{N+2}}
+\gm_o\bigg(\frac{\theta}{\mu^-+a\xi\om}\bigg)^{\frac2{N+2}}\bigg].
\]
This finishes the proof.\hfill\bbox\\
\\
Some remarks are in order.
\begin{remark}
Without loss of generality, we may assume that $\mu^-<\frac12\xi\om$. In such a case
the quantity $A$ above reduces to
\[
A=\frac{\gm}{(1-a)^2}\bigg[\bigg(\frac{\xi\om}{\theta}\bigg)^{\frac{N}{N+2}}
+\bigg(\frac{\theta}{a\xi\om}\bigg)^{\frac2{N+2}}\bigg].
\]
\end{remark}
\begin{remark}
The either-or conclusion is necessary. Without $\xi\om>I_{p,\rho}(y,s)$, in general
 one cannot obtain
 \[
 u\ge \mu^-+a\xi\om\quad\text{a.e. in }(y,s)+Q_{\frac\rho2}(\theta).
 \]
 See Remark C.1 in Appendix C of \cite{L}.
\end{remark}
\begin{lemma}\label{Lm:4:2}
Let $u$ be a non-negative, 
locally bounded, local, weak sub-solution to equation 
(\ref{Eq:1:1}), in $E_T$. Assume that
\begin{equation}\label{Eq:B:6:1}
\om\ge{\frac{1}{b+1}}\mu_+,
\end{equation}
for some positive parameter $b$ to be chosen later.
There exists a positive number $\nu_+$, 
depending upon $\om$, $\theta$, $\xi$, $a$ and $N$, such that if
\begin{equation*}
\big|[u\ge \mu_+-\xi\om]\cap[(y,s)+Q_{2\rho}(\theta)]\big|\le
\nu_+|Q_{2\rho}(\theta)|
\end{equation*}
then 
\begin{equation*}
u\le \mu_+-a\xi\om\qquad\text{ a.e. in }\>
(y,s)+Q_\rho(\theta).
\end{equation*}
\end{lemma}
\noi{\bf Proof.} 
Assume $(y,s)=(0,0)$ and for $n=0,1,\dots$ set
\begin{equation*}
\rho_n=\rho+\frac{\rho}{2^n},\qquad K_n=K_{\rho_n},\qquad Q_n=K_n\times(-\theta\rho_n^2].
\end{equation*}
Let $\z$ be a non-negative, piecewise smooth cutoff function on 
$Q_n$ defined as in the previous lemma.
Introduce the sequence of truncating levels 
\begin{equation*}
k_n =\mu_+-\xi_n\om\qquad\text{ with }\qquad
\xi_n=a\xi+\frac{1-a}{2^n}\xi,
\end{equation*}
and write down the energy estimates (\ref{Eq:2:13}) over the 
cylinder $Q_n$, for the truncated function $\uknp$. Taking 
also into account (\ref{Eq:B:6:1}), this gives
\begin{equation*}
\begin{aligned}
\essup_{-\theta\rho_n^2<t\le 0}&\int_{K_n}\uknp^2\z^2(x,t)dx\\
&\quad + \iint_{Q_n}\frac{|D[\uknp\z]|^2}u dxd\tau\\
&\le\gm \frac{2^{2n}}{\rho^2}(\xi\om)^2
\iint_{Q_n} \Big(\frac{1}{(1-\xi)\om}
+\frac1{\theta}\Big)
\chi_{[u>k_n]}dxd\tau\\
&\le \gm \frac{2^{2n}}{\rho^2} (\xi\om)^2 
\frac1{(1-\xi)\om} 
\Big(1+\frac\om{\theta}\Big) |[u>k_n]\cap Q_n|.
\end{aligned}
\end{equation*}
To estimate below the second integral on the left-hand side, 
take into account that  $u\le\mu_+$ and (\ref{Eq:B:6:1}). This gives
\begin{equation*}
\begin{aligned}
\iint_{Q_n}\frac{|D[\uknp\z]|^2}u dxd\tau
\ge\,{\frac{1}{(b+1)\om}}
\iint_{Q_n}|D[\uknp\z]|^2 dxd\tau.
\end{aligned}
\end{equation*}
Setting
\begin{equation*}
A_n=[u>k_n]\cap Q_n \quad\text{ and }\quad Y_n=\frac{|A_n|}{|Q_n|},
\end{equation*}
and combining these estimates gives
\begin{equation}\label{Eq:B:6:2}
\begin{aligned}
\essup_{-\theta\rho_n^2<t\le  0}
\int_{K_n}&\uknp^2\z^2(x,t)dx\\
&\quad +\frac{1}{(b+1)\om}\iint_{Q_n}|D[\uknp\z]|^2 
dxd\tau\\
&\le \gm\frac{2^{2n}}{\rho^2} 
\frac{(\xi\om)^2}{(1-\xi)\om} 
\Big(1+\frac{\om}{\theta}\Big)\,|A_n|.
\end{aligned}
\end{equation}
Apply H\"older's inequality and the embedding 
Proposition~3.1 of Chapter 1 of \cite{DB}, 
and 
recall that $\z=1$ on $Q_{n+1}$, to get 
\begin{equation*}
\begin{aligned}
\Big(\frac{1-a}{2^{n+1}}\Big)^2(\xi\om)^2&|A_{n+1}|
\le \iint_{Q_{n+1}}\uknp^2dxd\tau\\
&\le\left(\iint_{Q_{n}}[\uknp\z]^{2\frac{N+2}{N}}
dxd\tau\right)^{\frac{N}{N+2}}
|A_n|^{\frac{2}{N+2}}\\
&\le \gm\Big(\iint_{Q_{n}}|D[\uknp\z]|^2dxd\tau
\Big)^{\frac{N}{N+2}}\\
&\quad\times\Big(\essup_{-\theta\rho_n^2<t\le 
0}\int_{K_{n}}[\uknp\z]^2(x,t)dx\Big)^{\frac{2}{N+2}}
|A_n|^{\frac{2}{N+2}}
\end{aligned}
\end{equation*}
for a constant $\gm$ depending only upon $N$. 
Combine this with (\ref{Eq:B:6:2}) to get
\begin{equation*}
\begin{aligned}
|A_{n+1}|\le&\frac{\gm\,2^{4n}}{(1-a)^2\rho^2}
\frac{(b+1)^{\frac N{N+2}}}{(1-\xi)}\frac1{\om^{\frac2{N+2}}}
\Big(1+\frac\om{\theta}\Big)|A_n|^{1+\frac{2}{N+2}}.
\end{aligned}
\end{equation*}
In terms of $Y_n=\frac{|A_n|}{|Q_n|}$ this can be rewritten as
\begin{eqnarray*}
Y_{n+1}&\le &\frac{\gm 2^{4n}(b+1)^{\frac N{N+2}}}{(1-a)^2(1-\xi)}
\Big(\frac\theta\om\Big)^{\frac{2}{N+2}}\Big(1+\frac\om\theta\Big)
\,Y_n^{1+\frac{2}{N+2}}.
\end{eqnarray*}
By Lemma~4.1 of Chapter 1 of \cite{DB}, 
$\{Y_n\}\to 0$ as $n\to\infty$, provided 
\begin{equation}\label{nu+}
\begin{aligned}
Y_o=\frac{|A_o|}{|Q_o|}\le& 
\left[\frac{(1-a)^2(1-\xi)}{\gm 4^{N+2}(b+1)^{\frac N{N+2}}}
\right]^{\frac{N+2}2}
\frac{\frac\om\theta}
{\big(1+\frac\om\theta\big)^\frac{N+2}{2}}\,\df{=}\,\nu_+.
\end{aligned}
\end{equation}
This finishes the proof.\hfill\bbox
\section{Proof of Theorem \ref{Thm:1}}
Fix $(x_o,t_o)\in E_T$ and let $\rho>0$ be so small that $(x_o,t_o)+Q_{\rho}\subset E_T$;
 we may assume that $(x_o,t_o)$ coincides with the origin. 
Set
\[\mu_+=\essup_{Q_{\rho}}u,\quad\mu_-=\essinf_{Q_{\rho}}u,\quad\om=\mu_+-\mu_-.\]
Without loss of generality we may assume $\om\le1$, such that
\[Q_{\rho}(\om)\subset Q_{\rho}\quad\text{and}\quad\essosc_{Q_{\rho}(\om)}u\le \om.\]
The proof now unfolds along several cases.
\subsection{Case I}
First of all, let us suppose 
\[\mu_-\ge\frac18\om\quad\Leftrightarrow\quad \mu_+\le9\mu_-.\]
Without loss of generality we may assume $\mu_+\le1$ such that
\[
Q_{\rho}(\mu_+)\subset Q_{\rho}\quad\text{and}\quad\essosc_{Q_{\rho}(\mu_+)}u\le\om.
\]
Introduce the change of time variable and unknown function
\[
\tau=\mu_+^{-1}t\quad\text{and}\quad v(\cdot,\tau)=\frac{u(\cdot,t)}{\mu_+}.
\]
Then 
\[
v_\tau-\dvg\frac{Dv}v=0\quad\text{weakly in}\quad Q_{\rho}
\]
with
\[\frac19\le v\le 1.\]
Thus, by the classical parabolic theory (\cite{LSU}), there exists $\eta\in(0,1)$
depending only on $N$ such that
\[\essosc_{Q_{\frac{\rho}2}}v\le(1-\eta)\essosc_{Q_{\rho}}v.\]
Returning to the original coordinates we conclude that
\[
\essosc_{Q_{\frac{\rho}2}(\om)}u\le\essosc_{Q_{\frac{\rho}2}(\mu_+)}u\le(1-\eta)\essosc_{Q_{\rho}(\mu_+)}u\le(1-\eta)\om.
\]
\subsection{Case II}
Now suppose
\[\mu_+>9\mu_-,\]
which is equivalent to 
\[\om>\frac89\mu_+.\]
Suppose in addition that 
\[|[u\le\mu^-+\frac12\om]\cap Q_{\rho}(\om)|\le \nu_-|Q_{\rho}(\om)|\]
where $\nu_-$ is defined in \eqref{mu} with $\xi=1/2$ and $\theta=\om$.
Then by Lemma \ref{Lm:4:1} with $a=\frac12$,
we have either
\[\om\le 2 I_{p,\rho}\]
or
\[u\ge\mu^-+\frac14\om\quad\text{ a.e. in }Q_{\frac{\rho}2}(\om).\]
The latter implies
\[\essosc_{Q_{\frac{\rho}2}(\om)}u\le\frac34\om.\]
\subsection{Case III}
As in the previous case suppose  that
\[\om>\frac89\mu_+,\]
but
\[
|[u\le\mu^-+\frac12\om]\cap Q_{\rho}(\om)|> \nu_-|Q_{\rho}(\om)|.
\]
Then there exists some 
\[-\om\rho^2\le s\le -\frac12\nu_-\om\rho^2\]
such that 
\begin{equation}\label{Eq:4:1}
|[u(\cdot,s)<\mu^-+\frac12\om]\cap K_{\rho}|>\frac12\nu_-|K_{\rho}|.
\end{equation}
Indeed, if the above inequality does not hold for any $s$ in the given 
interval, then 
\begin{align*}
|[u<\mu^-+\frac12\om]\cap Q_{\rho}(\om)|=
&\int_{-\om\rho^2}^{-\frac12\nu_-\om \rho^2}|[u(\cdot,s)<\mu^-+\frac12\om]\cap Q_{\rho}(\om)|\,ds\\
&+\int^{0}_{-\frac12\nu_-\om\rho^2}|[u(\cdot,s)<\mu^-+\frac12\om]\cap Q_{\rho}(\om)|\,ds\\
&\le\nu_- |Q_{\rho}(\om)|.
\end{align*}
Since $\mu_+-\frac14\om>\mu^-+\frac12\om$ always holds, \eqref{Eq:4:1} implies
\[|[u(\cdot,s)>\mu_+-\frac14\om]\cap K_{\rho}|\le(1-\frac12\nu_-)|K_{\rho}|.
\]
Based on this, we use the logarithmic estimate to show that such a measure
theoretical information propagates in time.

\noi{\bf Claim 1: } There exists a positive integer $n_*$ depending only on $N$ such that
\[|[u(\cdot,t)>\mu_+-\frac{\om}{2^{n_*}}]\cap K_{\rho}|\le(1-\frac14\nu_-^2)|K_{\rho}|\quad\text{ for all }s<t<0.\]
{\bf Proof.} In the logarithmic estimate we take 
\[
k=\mu_+-\frac14\om,\quad c=\frac{\om}{2^{n+2}}.
\]
This gives
\[
\psi(u)=\ln^+\bigg[\frac{H}{H-[u-(\mu_+-\frac14\om)_+]-\frac{\om}{2^{n+2}}}\bigg]
\]
where \[H=\essup_{K_{2\rho}\times(s,0)}(u-k)_+.\]
Choose a cutoff function $\z$ which satisfies
$\z=1$ on $K_{(1-\sig)\rho}$ and $\z=0$ on $\pl K_\rho$, 
such that \[|D\z|\le\frac1{\sig\rho}.\]
Hence, for all $s<t<0$,
\begin{align*}
\int_{K_{(1-\sig)\rho}(y)}
\psi^2(u)(x,t)dx
&\le\int_{K_\rho(y)}\psi^2(u)(x,s)dx+\frac{\gm}{(\sig\rho)^2} \int_s^0\int_{K_\rho}\frac{\psi(u)}u dxdt.
\end{align*}
Note that \[\psi\le n\ln2.\]
The first term on the right-hand side is estimated by
\[\int_{K_\rho(y)}\psi^2(u)(x,s)dx\le n^2\ln^2 2(1-\frac12\nu)|K_{\rho}|.\]
The second term is estimated by
\[\frac{\gm}{(\sig\rho)^2} \int_s^0\int_{K_\rho}\frac{\psi(u)}u dxdt
\le\frac{\gm n}{(\sig\rho)^2}(\om \rho^2)\om^{-1}|K_{\rho}|\le\frac{\gm n}{\sig^2}|K_{\rho}|.\]
The left-hand side is estimated below by integrating over the smaller set 
\[
[u(\cdot,t)>\mu_+-\frac{\om}{2^{n+2}}]\cap K_{(1-\sig)\rho}.
\]
On such a set \[\psi^2\ge\ln^2\bigg(\frac{\frac{\om}4}{\frac{\om}{2^{n+1}}}\bigg)=(n-1)^2\ln^2 2.\]

Thus, combining all above estimates yields
\[|[u(\cdot,t)>\mu_+-\frac{\om}{2^{n+2}}]\cap K_{(1-\sig)\rho}|
\le\bigg(\frac{n}{n-1}\bigg)^2(1-\frac12\nu_-)|K_\rho|+\frac{\gm}{n\sig^2}|K_{\rho}|\]
for all $s<t<0$. On the other hand,
\begin{align*}
&|[u(\cdot,t)>\mu_+-\frac{\om}{2^{n+2}}]\cap K_\rho|\\
&\le|[u(\cdot,t)>\mu_+-\frac{\om}{2^{n+2}}]\cap K_{(1-\sig)\rho}|+|K_\rho-K_{(1-\sig)\rho}|\\
&\le|[u(\cdot,t)>\mu_+-\frac{\om}{2^{n+2}}]\cap K_{(1-\sig)\rho}|+N\sig |K_\rho|.
\end{align*}
Therefore,
\[
|[u(\cdot,t)>\mu_+-\frac{\om}{2^{n+2}}]\cap K_\rho|
\le\bigg[\bigg(\frac{n}{n-1}\bigg)^2(1-\frac12\nu_-)+\frac{\gm}{n\sig^2}+N\sig\bigg]|K_\rho|.
\]
The claim is proved by choosing $\sig$ so small, and then $n$ large enough.\hfill\bbox
\\
\\

Using the measure theoretical information obtained for every time level of the cylinder
\[Q_\rho(\frac12\nu_-\om)=K_{\rho}\times(-\frac12\nu_-\om\rho^2,0]\] in {\bf Claim 1}, we
are able to show

\noi{\it {\bf Claim 2:} For any $\nu_*\in(0,1)$ there exists a positive integer $l$ such that
\[|[u>\mu_+-\frac{\om}{2^{n_*+l}}]\cap Q_\rho(\frac12\nu_-\om)|\le\nu_*|Q_\rho(\frac12\nu_-\om)|.\]}
{\bf Proof.} Let $Q=Q_\rho(\frac12\nu_-\om)$ and $Q'=Q_{2\rho}(\frac12\nu_-\om)$.
Write down the energy estimate over $Q'$ for 
\[(u-k_j)_+\quad\text{where }k_j=\mu_+-\frac{\om}{2^j}\quad\text{for }j=n_*,\,\cdots,\,n_*+l.\]
Choose a cutoff function $\z$ which satisfies $\z=1$ on $Q$, and 
vanishes on the parabolic boundary of $Q'$, such that
\[
|D\z|\le\frac1\rho,\quad |\z_t|\le\frac1{\nu_-\om\rho^2}.
\]
Then keeping in mind  we assumed at the beginning that $\om>8\mu_+/9$, the energy estimate gives
\begin{equation}\label{graduk}
\om^{-1}\iint_Q|D(u-k_j)_+|^2\,dxdt\le\frac{\gm}{\nu_-\om\rho^2}\bigg(\frac{\om}{2^j}\bigg)^2|Q|.
\end{equation}
Next, we apply the discrete isoperimetric inequality on page 15 of \cite{DB}
to $u(\cdot,t)$ for $-\frac12\nu_-\om\rho^2<t<0$, over the cube $K_\rho$,
for levels $k_j<k_{j+1}$.
Taking into account the measure theoretical information from Claim 1,
this gives
\begin{align*}
\frac{\om}{2^{j+1}}&|[u(\cdot,t)>k_{j+1}]\cap K_\rho|\\
&\le\frac{\gm \rho^{N+1}}{|[u(\cdot,t)<k_j]\cap K_\rho|}\int_{[k_j<u(\cdot,t)<k_{j+1}]\cap K_{\rho}}|Du|\,dx\\
&\le\frac{\gm}{\nu_-^2}\rho\bigg(\int_{[k_j<u(\cdot,t)<k_{j+1}]\cap K_{\rho}}|Du|^2\,dx\bigg)^{\frac12}\\
&\quad\times|([u(\cdot,t)>k_j]-[u(\cdot,t)>k_{j+1}])\cap K_\rho|^\frac12.
\end{align*}
Set
\[
A_j=[u>k_j]\cap Q
\]
and integrate the above estimate in $dt$ over $-\frac12\nu_-\om\rho^2<t<0$;
we obtain
\[
\frac{\om}{2^j}|A_{j+1}|\le\frac{\gm}{\nu_-^2}\rho\bigg(\iint_Q|D(u-k_j)_+|^2\,dxdt\bigg)^\frac12(|A_j|-|A_{j+1}|)^\frac12.
\]
Now square both sides of the above inequality and use \eqref{graduk}
to estimate the term containing $D(u-k_j)_+$, to obtain
\[|A_{j+1}|^2\le\frac{\gm}{\nu_-^5}|Q|(|A_j|-|A_{j+1}|).\]
Add these inequalities from $n_*$ to $n_*+l-1$ to obtain
\[
(l-2)|A_{n_*+l}|^2\le\sum_{j=n_*+1}^{n_*+l-1}|A_{j+1}|^2\le\frac{\gm}{\nu_-^5}|Q|^2.
\]
From this 
\[|A_{n_*+l}|\le\frac1{\sqrt{l-2}}\sqrt{\frac{\gm}{\nu_-^5}}|Q|.\]
Given a number $\nu_*\in(0,1)$, we can choose $l$ large enough to guarantee
\[\frac1{\sqrt{l-2}}\sqrt{\frac{\gm}{\nu_-^5}}<\nu_*.\]
This finishes the proof.\hfill\bbox

Now we are ready to finish {\bf Case III}. 
Choosing $a=\frac12$, $b=\frac18$ and $\theta=\frac12\nu_-\om$, the constant $\nu_+$
from \eqref{nu+} becomes
\[
\nu_+=\frac1{\gm}\nu_-^{\frac{N}2}(1-\xi)^{\frac{N+2}2}
\]
with 
\[
\xi=\frac1{2^{n_*+l}},\quad\text{ $n_*$ and $l$ to be determined.}
\]
 Then after choosing $n_*$ from {\bf Claim 1}, we can choose $l$ so large that
 \[
 \frac1{\sqrt{l-2}}\sqrt{\frac{\gm}{\nu_-^5}}<\frac1{\gm}\nu_-^{\frac{N}2}\bigg(1-\frac1{2^{n_*+l}}\bigg)^{\frac{N+2}2}.
 \]
 By Lemma \ref{Lm:4:2}, we obtain that
\[u\le \mu_+-\frac{\om}{2^{n_*+l+1}}\quad\text{in}\quad Q_{\frac\rho2(\frac12\nu_-\om)}.\]
This in turn implies
\[
\essosc_{Q_{\frac\rho2(\frac12\nu_-\om)}}u\le (1-\dl)\om\quad\text{where}\quad\dl=\frac{1}{2^{n_*+l+1}}.
\]
{\bf Proof of Theorem \ref{Thm:1}.}
Combining all these cases above, we have proved that once we have
$$\essosc_{Q_{\rho}(\om)}u\le\om,$$
 we can find  positive constants 
$$c=\min\{\frac12,\,\sqrt{\frac{\nu_-}8}\}\quad\text{and}\quad \lm=\min\{\frac34,\,1-\eta,\,1-\dl\}$$ such that
\[
\essosc_{Q_{c\rho}(\lm\om)}u\le \max\{\lm\om,\, 2I_{p,\rho}\}.
\]
Relabel the quantities $\rho$ and $\om$ chosen above as $\rho_o$ and $\om_o$. 
Now let $$\om_1=\max\{\lm\om_o,\, 2I_{p,\rho_o}\}\quad\text{and}\quad\rho_1=c\rho_o$$
 such that
\[
Q_{\rho_1}(\om_1)\subset Q_{c\rho_o}(\lm\om_o)\quad\text{and}\quad 
\essosc_{Q_{\rho_1}(\om_1)}u\le\om_1.
\]
The above set inclusion is verified if 
\[
\om_1\rho^2_1=\max\{\lm\om_o,\, 2I_{q,\rho_o}\}(c\rho_o)^2\le \lm\om_o(c\rho_o)^2.
\]
This holds naturally unless
\[
\lm\om_o\le 2I_{p,\rho_o},
\]
but then there is nothing to prove since
\[
\essosc_{Q_{\rho_o}(\om_o)}u\le \frac{2}{\lm}I_{p,\rho_o}.
\]
Hence, with such a choice of $c$ we now have 
\[
\essosc_{Q_{\rho_1}(\om_1)}u\le\om_1.
\]
According to what we have shown, one has
\[
\essosc_{Q_{c\rho_1}(\lm\om_1)}u\le \max\{\lm\om_1,\, 2I_{p,\rho_1}\}.
\]
Now define
\[
\om_2=\max\{\lm\om_1,\, 2I_{p,\rho_1}\}\quad\text{and}\quad\rho_2=c\rho_1,
\]
and we want to show that
\[
Q_{\rho_2}(\om_2)\subset Q_{c\rho_1}(\lm\om_1)\quad\text{and hence}\quad 
\essosc_{Q_{\rho_2}(\om_2)}u\le\om_2.
\]
The above set inclusion is verified if 
\[
\om_2\rho^2_2=\max\{\lm\om_1,\, 2I_{p,\rho_1}\}(c\rho_1)^2\le \lm\om_1(c\rho_1)^2.
\]
This holds naturally unless
\[
\lm\om_1\le 2I_{p,\rho_1},
\]
but then there is nothing to prove since
\[
\essosc_{Q_{\rho_1}(\om_1)}u\le \frac{2}{\lm}I_{p,\rho_1}.
\]

Iterating in this fashion, one concludes that there are positive numbers
$c,\,\lm\in(0,1)$ such that constructing
\[
\rho_n=c^n\rho_o \quad\text{and}\quad\om_n=\max\{\lm\om_{n-1},\, 2I_{p,\rho_{n-1}}\}
\]
one obtains
\[
\essosc_{Q_{\rho_n}(\om_n)}u\le\max \{\om_n,\,\frac{2}{\lm}I_{p,\rho_{n-1}}\}
\]
Let $0<r<\rho_o\le R_o$ be fixed.
Since the sequence $\{\om_n\rho_n^2\}$
is strictly decreasing and gives a partition of the interval
$(0,\,\om_o\rho_o^2)$, there must be some positive integer $n$ such that
\[
\rho_{n+1}^2\om_{n+1}\le r^2\om_o<\rho_n^2\om_n.
\]
Noting that 
$$\om_n\le\om_o\quad\text{and}\quad \om_{n+1}\ge\lm^{n+1}\om_o,$$ 
this implies that $r<\rho_n$ and 
\[
n+1\ge\frac{\ln\big(\frac{r}{\rho_o}\big)}{\ln(\sqrt{\lm}c)}.
\]
Then it is not hard to see that
\[
Q_{r}(\om_o)\subset Q_{\rho_n}(\om_n),
\]
and either
\[
\essosc_{Q_{r}(\om_o)}u\le \om_n\quad\text{or}\quad \om_{n-1}\le\frac{2}{\lm}I_{p,\rho_{n-1}}.
\]
Note that
\begin{align*}
\om_n&=\max\{\lm^n\om_o,2\lm^{n-1}I_{p,\rho_o},\cdots,2\lm I_{p,\rho_{n-2}},2I_{p,\rho_{n-1}}\}\\
&\le \lm^n\om_o+\frac2{1-\lm}I_{p,\rho_o}.
\end{align*}
Here we made the convention that  the function $\rho\to I_{p,\rho}$ is non-decreasing.
Otherwise, one could use
\[
\tilde{I}_{p,\rho}=\sup_{0<\tau<\rho}I_{p,\tau}.
\]
Thus, there is some $\al\in(0,1)$ depending only on the data such that
\[
\essosc_{Q_{r}(\om_o)}u\le\lm^n\om_o+C I_{p,\rho_o}\le 
\bar{C}\bigg[\om_o\bigg(\frac{r}{\rho_o}\bigg)^\al+I_{p,\rho_o}\bigg]
\]
where
\[
\al=\frac{\ln\lm}{\ln(\sqrt{\lm}c)}\quad\text{and}\quad C=\max\left\{\frac2{1-\lm},\frac{2}{\lm}\right\}.
\]
Now choose $\rho_o=R_o^{1-\mu}r^\mu$, and conclude we have
\[
\essosc_{Q_{r}(\om_o)}u\le\bar{C}\bigg[\om_o\bigg(\frac{r}{R_o}\bigg)^{(1-\mu)\al}+I_{p,R^{1-\mu}_o r^\mu}\bigg].
\]

Naian Liao\\
College of Mathematics and Statistics\\
Chongqing University\\
Chongqing, China, 401331\\
email: {\tt liaon@cqu.edu.cn}
\end{document}